\documentclass[12pt,thmsa]{article}
\usepackage{amsfonts}

\input{tcilatex}
\begin{document}

\author{Harald Grosse$^a$ \and Karl-Georg Schlesinger$^b$ \qquad \\
$^a$Institute for Theoretical Physics\\
University of Vienna\\
Boltzmanngasse 5\\
A-1090 Vienna, Austria\\
e-mail: grosse@doppler.thp.univie.ac.at\\
$^b$Erwin Schr\"{o}dinger Institute for Mathematical Physics\\
Boltzmanngasse 9\\
A-1090 Vienna, Austria\\
e-mail: kgschles@esi.ac.at}
\title{On a noncommutative deformation of the Connes-Kreimer algebra}
\date{}
\maketitle

\begin{abstract}
We study a noncommutative deformation of the commutative Hopf algebra $%
\mathcal{H}_R$ of rooted trees which was shown by Connes and Kreimer to
describe the mathematical structure of renormalization in quantum field
theories. The requirement of the existence of an antipode for the
noncommutative deformation leads to a natural extension of the algebra.
Noncommutative deformations of $\mathcal{H}_R$ might be relevant for
renormalization of field theories on noncommutative spaces and there are
indications that in this case the extension of the algebra should be linked
to a mixing of infrared and ultraviolet divergences. We give also an
argument that for a certain class of noncommutative quantum field theories
renormalization should be linked to a noncommutative and noncocommutative
self-dual Hopf algebra which can be seen as a noncommutative counterpart of
the Grothendieck-Teichm\"{u}ller group.
\end{abstract}

\section{Introduction}

The Hopf algebra $\mathcal{H}_R$ of rooted trees (see \cite{CK}) can be seen
as the abstract mathematical structure behind the renormalization scheme
employed by physicists in quantum field theory (see \cite{CK} and the
literature cited there). The study of quantum field theories on
noncommutative spaces has started in recent years (see \cite{CDP 1998}, \cite
{CDP 2000}, \cite{GMS}, \cite{MSSW}, \cite{Oec}) but the question of
renormalization of such theories is still an open problem. Since the Hopf
algebra $\mathcal{H}_R$ - by the nature of an algebra of polynomials - is
commutative, it is natural to ask if noncommutative deformations of $%
\mathcal{H}_R$ could be the proper algebraic setting for renormalization of
field theories on noncommutative spaces. This is the motivation for studying
a very simple noncommutative deformation of $\mathcal{H}_R$, here. If this
deformation would indeed be linked to renormalization on noncommutative
spaces, we find indications that a mixing of infrared and ultraviolet
divergences should occur, as has been observed in examples (see \cite{MRS}).
We focus, here, on the algebraic properties of the deformation and give only
a brief discussion of a possible realization of the deformed algebra in
terms of $q$-integrals and the shift and particle number operators.

Besides the question of renormalization, a second motivation for the study
of noncommutative deformations of $\mathcal{H}_R$ comes from considerations
on trialgebraic deformations of Hopf algebras and a noncommutative and
noncocommutative Hopf algebra $\mathcal{H}_{GT}$ taking the role which the
Grothendieck-Teichm\"{u}ller group plays for quasitensor categories, there
(see \cite{Sch}). In section 3 we give an argument from the algebraic
properties of noncommuative quantum field theories that for a certain class
of such theories renormalization should, indeed, be expected to be linked to 
$\mathcal{H}_{GT}$.

\bigskip

\section{The deformation of $\mathcal{H}_R$}

For a rooted tree $t$ as defined in \cite{CK}, let $v(t)$ be the number of
vertices of $t$. We make the convention that for the unit element $e$ of $%
\mathcal{H}_R$ 
\[
v\left( e\right) =0 
\]
For a monomial of rooted trees, the number of vertices $v$ is, of course,
defined as the sum of the vertex numbers of its factors. Let $q\neq 1$ be a
complex deformation parameter. We introduce a deformation $\mathcal{H}_{R,q}$
of $\mathcal{H}_R$ by replacing commutativity, i.e. 
\[
t_1t_2=t_2t_1 
\]
for rooted trees by the condition 
\begin{equation}
t_1t_2=q^{v\left( t_2\right) -v\left( t_1\right) }t_2t_1  \label{1}
\end{equation}
In addition, we require that $e$ is invertible, i.e. there exists $e^{-1}$
with 
\[
e^{-1}e\ t=e\ e^{-1}t=t\ e^{-1}e=t\ e\ e^{-1}=t 
\]
for all rooted trees $t$. Clearly, this defines a deformation of the
associative algebra structure of $\mathcal{H}_R$. We have to study now how
this deformation effects the other structural elements involved in $\mathcal{%
H}_R$. Observe, first, that (\ref{1}) involves a deformation of the unit
element $e$, too, since 
\begin{equation}
e\ t=q^{v\left( t\right) }t\ e  \label{2}
\end{equation}
for an arbitrary rooted tree $t$. Let us take a look at the coproduct $%
\Delta $, next. Remember that $\Delta $ can be defined for a rooted tree $t$
as (see\cite{CK}) 
\begin{equation}
\Delta \left( t\right) =e\otimes t+t\otimes e+\sum_CP^C\left( t\right)
\otimes R^C\left( t\right)  \label{3}
\end{equation}
where the sum is taken over admissible cuts $C$ of $t$ and $R^C\left(
t\right) $ is the rooted tree containing the root of $t$ while $P^C\left(
t\right) $ denotes the complementary monomial of rooted trees.

\bigskip

\begin{lemma}
$\Delta $ defines a coproduct also for the deformed algebra as given by (\ref
{1}).
\end{lemma}

\proof%
We have to check that (\ref{1}) is consistent with $\Delta $, i.e. we have
to check that 
\[
\Delta \left( t_1\right) \Delta \left( t_2\right) =q^{v\left( t_2\right)
-v\left( t_1\right) }\Delta \left( t_2\right) \Delta \left( t_1\right) 
\]
for rooted trees $t_1$, $t_2$.

From (\ref{3}) it follows that 
\begin{eqnarray*}
&&\Delta \left( t_1\right) \Delta \left( t_2\right) \\
&=&\left( e\otimes t_1\right) \left( e\otimes t_2\right) +\left( e\otimes
t_1\right) \left( t_2\otimes e\right) +\left( e\otimes t_1\right) \left(
\sum_CP^C\left( t_2\right) \otimes R^C\left( t_2\right) \right) \\
&&+\left( t_1\otimes e\right) \left( e\otimes t_2\right) +\left( t_1\otimes
e\right) \left( t_2\otimes e\right) +\left( t_1\otimes e\right) \left(
\sum_CP^C\left( t_2\right) \otimes R^C\left( t_2\right) \right) \\
&&+\left( \sum_CP^C\left( t_1\right) \otimes R^C\left( t_1\right) \right)
\left( e\otimes t_2\right) +\left( \sum_CP^C\left( t_1\right) \otimes
R^C\left( t_1\right) \right) \left( t_2\otimes e\right) \\
&&+\left( \sum_CP^C\left( t_1\right) \otimes R^C\left( t_1\right) \right)
\left( \sum_CP^C\left( t_2\right) \otimes R^C\left( t_2\right) \right)
\end{eqnarray*}
Using the fact that 
\[
v\left( t\right) =v\left( P^C\left( t\right) \right) +v\left( R^C\left(
t\right) \right) 
\]
for any rooted tree $t$ and equation (\ref{2}), the desired result follows.%
\endproof%

\bigskip

The counit $\varepsilon $ is defined as 
\[
\varepsilon \left( t\right) =0 
\]
for $t\neq e$ and 
\[
\varepsilon \left( e\right) =1 
\]
in\cite{CK}.

\bigskip

\begin{lemma}
The counit $\varepsilon $ is compatible with (\ref{1}).
\end{lemma}

\proof%
By calculation.%
\endproof%

\bigskip

\begin{corollary}
$\mathcal{H}_{R,q}$ has the structure of a nonunital bialgebra.%
\endproof%
\end{corollary}

\begin{remark}
By results of \cite{GeSch}, there exists a deformation of $\mathcal{H}_R$
equivalent to $\mathcal{H}_{R,q}$ which leaves $e$ fixed, too, i.e. there
exists an equivalent deformation of $\mathcal{H}_R$ into a full bialgebra
(here, equivalence of deformations refers to the deformation of the
associative product and the coassociative coproduct).
\end{remark}

\bigskip

Finally, let us consider the question of the existence of an antipode.
Suppose, $S$ would be an antipode for $\mathcal{H}_{R,q}$. Since an antipode
is always an algebra antihomomorphism (see e.g. \cite{KS}), it follows from (%
\ref{2}) that 
\begin{equation}
S\left( t\right) S\left( e\right) =q^{v\left( t\right) }S\left( e\right)
S\left( t\right)  \label{4}
\end{equation}
for any rooted tree $t$. Suppose now that $q\neq -1$. By the definition of $%
\mathcal{H}_{R,q}$, this implies that $S(e)$ would have to involve monomials
of rooted trees with a larger number of vertices than any fixed monomial of
rooted trees which is, obviously, a contradiction. So, there can not exist
an antipode on $\mathcal{H}_{R,q}$ for $q\neq -1$.

\bigskip

\begin{remark}
For $q=-1$ no contradiction does arise, here, because $q=q^{-1}$, then. It
is interesting to observe that the algebraic deformation theory of $\mathcal{%
H}_R$ immediately seems to mirror the fact that the renormalization scheme,
described by $\mathcal{H}_R$, can be generalized without problems to the
supersymmetric setting.
\end{remark}

\bigskip

Equation (\ref{4}) also points a way to a partial solution of the problem
how to introduce an antipode in the deformed case. In a certain sense, we
simply have to consider $\mathcal{H}_{R,q}$ and $\mathcal{H}_{R,q^{-1}}$ at
once. Denote rooted trees in $\mathcal{H}_{R,q^{-1}}$ by $\widehat{t}$ to
distinguish them from those of $\mathcal{H}_{R,q}$, and let $\widetilde{%
\mathcal{H}_{R,q}}$ be the nonunital bialgebra generated from 
\[
\mathcal{H}_{R,q}\oplus \mathcal{H}_{R,q^{-1}} 
\]
with the relations 
\begin{equation}
q^{v\left( t_1\right) }t_1\widehat{t_2}=\widehat{t_1t_2}=q^{-v\left(
t_2\right) }\widehat{t_1}t_2  \label{5}
\end{equation}
imposed. Equation (\ref{5}) assures that the counit and coproduct defined
separately for $\mathcal{H}_{R,q}$ and $\mathcal{H}_{R,q^{-1}}$ can be
consistently combined into a single counit and coassociative coproduct for $%
\widetilde{\mathcal{H}_{R,q}}$.

Denote now by $S$ the antipode of $\mathcal{H}_R$ which is defined as (see 
\cite{CK}): 
\[
S\left( e\right) =e 
\]
and 
\begin{equation}
S\left( t\right) =-t-\sum_CS\left( P^C\left( t\right) \right)
e^{-1}R^C\left( t\right)  \label{6}
\end{equation}
Observe that we have inserted the element $e^{-1}$ on the right hand side of
(\ref{6}). This is not necessary for the case of $\mathcal{H}_R$, where $e$
is the unit, but we have to keep track of this appearance of $e^{-1}$ in the
noncommutative case.

Define $S_q$ on $\widetilde{\mathcal{H}_{R,q}}$ by 
\[
S_q\left( t\right) =\widehat{S\left( t\right) } 
\]
and 
\[
S_q\left( \widehat{t}\right) =S\left( t\right) 
\]
i.e. $S_q$ leads to an interchange of $\mathcal{H}_{R,q}$ and $\mathcal{H}%
_{R,q^{-1}}.$

\bigskip

\begin{lemma}
$S_q$ defines a left antipode for $\widetilde{\mathcal{H}_{R,q}}$.
\end{lemma}

\proof%
We have (denoting the multiplication by the product by $m$) 
\begin{eqnarray*}
&&m\left[ \left( S_q\otimes id\right) \Delta \left( t\right) \right] \\
&=&\widehat{e}\ t-\widehat{t}\ e-\sum_C\widehat{S}\left( P^C\left( t\right)
\right) \widehat{e^{-1}}\widehat{R^C\left( t\right) \ }e+\sum_C\widehat{S}%
\left( P^C\left( t\right) \right) R^C\left( t\right) \\
&=&q^{v\left( t\right) }\widehat{e}\ \widehat{t}-\widehat{t}\ \widehat{e}%
-\sum_C\widehat{S}\left( P^C\left( t\right) \right) \widehat{e^{-1}}\widehat{%
R^C\left( t\right) \ }\widehat{e}+\sum_Cq^{v\left( R^C\left( t\right)
\right) }\widehat{S}\left( P^C\left( t\right) \right) \widehat{R^C\left(
t\right) } \\
&=&0=\varepsilon \left( t\right)
\end{eqnarray*}
where we have used the definition of $e^{-1}$ and equation (\ref{5}). 
\endproof%

\bigskip

On the other hand, we have 
\begin{eqnarray*}
&&m\left[ \left( id\otimes S_q\right) \Delta \left( t\right) \right] \\
&=&t\ \widehat{e}-e\ \widehat{t}-\sum_Ce\ \widehat{S}\left( P^C\left(
t\right) \right) \widehat{e^{-1}}\widehat{R^C\left( t\right) }%
+\sum_CP^C\left( t\right) \widehat{S}\left( R^C\left( t\right) \right) \\
&=&\widehat{e}\ \widehat{t}-\widehat{e}\ \widehat{t}-\sum_Cq^{-v\left(
P^C\left( t\right) \right) }\widehat{S}\left( P^C\left( t\right) \right) \ 
\widehat{e}\ \widehat{e^{-1}}\widehat{R^C\left( t\right) }+\sum_Cq^{-v\left(
P^C\left( t\right) \right) }\widehat{P^C\left( t\right) }\widehat{S}\left(
R^C\left( t\right) \right) \\
&=&\sum_Cq^{-v\left( P^C\left( t\right) \right) }\left( \widehat{P^C\left(
t\right) }\widehat{S}\left( R^C\left( t\right) \right) -\widehat{S}\left(
P^C\left( t\right) \right) \widehat{R^C\left( t\right) }\right)
\end{eqnarray*}
which is non vanishing for $q\neq 1$ . So, $S_q$ does not define a right
antipode and, therefore, not an antipode for $\widetilde{\mathcal{H}_{R,q}}$.

On the other hand, one has the following result:

\bigskip

\begin{lemma}
There is a deformation of $\mathcal{H}_R$ equivalent (in the sense of
deformations of the associative product and coassociative coproduct) to $%
\mathcal{H}_{R,q}$ which carries a full Hopf algebra structure.
\end{lemma}

\proof%
As we remarked already above, there is a deformation equivalent to $\mathcal{%
H}_{R,q}$ which is a unital bialgebra. But then - since $\mathcal{H}_R$ is a
Hopf algebra - an antipode exists for this deformation (see \cite{CP}).%
\endproof%

\bigskip

It is, in general, not a straightforward task to explicitly construct the
equivalent deformations given in existence results of the kind above. So,
having an explicitly given left antipode, one could try to work some way
with this. The most prominent feature of $S_q$ is the interchange of $%
\mathcal{H}_{R,q}$ and $\mathcal{H}_{R,q^{-1}}$. In the renormalization
scheme the antipode corresponds to the construction of counter terms for
subdivergences of integrals. One notes that for the $q$-integrals as defined
in the setting of $q$-calculus (which often gives a good toy model for the
truely noncommutative situation) the exchange 
\[
q\longleftrightarrow q^{-1} 
\]
leads - modulo a constant term and a constant scaling factor - to a formal
exchange of $q$-integrals 
\[
\int_c^\infty d_qf\longleftrightarrow \int_0^cd_{q^{-1}}f 
\]
(see e.g. \cite{KS} for the definition of the $q$-integral), i.e. if the
left antipode introduced above is relevant for renormalization on
noncommutative spaces, this might be a hint that a renormalization scheme
based on $S_q$ should involve a mixing of ultraviolet and infrared
divergences.

This suggests that we can - at least in a formal sense - realize the algebra 
$\widetilde{\mathcal{H}_{R,q}}$ as follows: For an element of the component $%
\mathcal{H}_{R,q}$, we write the corresponding integral of $\mathcal{H}_R$
but understand all integrals to be replaced by $q$-integrals and we assume
that the formal variable of integration of the outer integral $q$-commutes
with the variables of integration of all the subintegrals (which all commute
with each other), i.e. for an $n$-fold integral we write 
\[
\int dy\ dx_1...dx_{n-1} 
\]
where all the $dx_i$ commute with each other and for $i=1,...,n-1$ 
\[
dy\ dx_i=q\ dx_i\ dy 
\]
It is easy to check that 
\[
\int dy\ dx_1...dx_{n-1}\ \int dy\ dx_1...dx_{m-1}=q^{m-n}\ \int dy\
dx_1...dx_{m-1}\ \int dy\ dx_1...dx_{n-1} 
\]
which is just the exchange rule required by the algebra $\mathcal{H}_{R,q}$.
For the elements of $\mathcal{H}_{R,q^{-1}}$, we use the corresponding
integrals with $q$ replaced by $q^{-1}$ (i.e. using integrals $\int_0^c...$
instead of $\int_c^\infty ...$). The use of the $q$-integral gives the
correct exchange between the components $\mathcal{H}_{R,q}$ and $\mathcal{H}%
_{R,q^{-1}}$ and since the $q$-integral is a kind of discrete approximation
to the usual integral, the qualitative structure of the singularities should
not change. By the nature of $S_q$, counter terms always alternate between
the ultraviolet and the infrared case. More concretely, we put the same
function into the integrand as in the toy model in \cite{CK} (i.e. $\frac
1{x+c}$ for the single integral, $\frac 1{\left( x_1+c\right) \left(
x_1+x_2\right) }$ for the double one, etc.) but replace the coordinate $x_i$
by $\frac 1{x_i}$ in every second $q$-integral. So, this should be a formal
realization of $\widetilde{\mathcal{H}_{R,q}}$.

Using the representation for the coordinates $x,y$ of the Manin plane on an
infinite dimensional space with basis $\left\{ \left| n\right\rangle ,n\in 
\Bbb{N}\right\} $, by 
\[
y\left| n\right\rangle =\left| n+1\right\rangle 
\]
and 
\[
x\left| n\right\rangle =q^n\left| n\right\rangle 
\]
we get the following interpretation of the above integrals: Each subintegral
involves - besides the $q$-integration - an application of the shift
operator which is quite natural if we view the rooted trees as quantum
objects, now, because with moving up the rooted tree the number of vertices,
taken into consideration in the integration process so far, increases by one
(i.e. which the shift recognizes by increasing the ``particle number'' by
one). For the outermost integral, we take the $q$-integration, again, but
then simply apply the particle number operator in an exponentiated form
(there is no more possibility to shift and we simply ``measure'' the final
particle number, now).

It would certainly be interesting to see a more detailed physical
realization of such a model.

\bigskip 

\begin{remark}
In terms of the generators $\delta _n$ (see \cite{CK}) the deformation of
the Connes-Kreimer algebra we have given above can be written as follows:
\[
\delta _n\delta _m=q^{m-n}\ \delta _m\delta _n
\]
With
\[
\delta _n=xy^n
\]
and $x,y$ as above, we get a representation of the deformed Connes-Kreimer
algebra in terms of the generators.
\end{remark}

\bigskip

\section{Bialgebra categories and noncommutative QFT}

In this section, we discuss some general algebraic properties of
noncommutative quantum field theories (i.e. quantum field theories on
noncommutative spaces, henceforth ncQFTs, for short) and give an abstract
argument why a noncommutative and noncocommutative Hopf algebra - albeit a
much more complicated one than the simple toy model $\widetilde{\mathcal{H}%
_{R,q}}$ - should be linked to renormalization of such theories.

Formally, the fields of a ncQFT can be seen as functions $\Phi $ with values 
\[
\Phi \left( \widehat{t},\widehat{x}\right) 
\]
in a noncommutative algebra $B$ and the variables in another noncommutative
algebra $A$ (see e.g. \cite{CDP 1998}). In general, these functions will not
be linear, especially, they are not restricted to the class of algebra
morphisms.

Suppose now that both algebras satisfy the Poincare-Birkhoff-Witt property
and can therefore be seen as arising from star-products. By \cite{CF}, \cite
{Kon}, both algebras can then be seen as arising from two dimensional
conformal field theories, i.e. - using the algebraic formulation of low
dimensional QFTs - we have two quasitensor categories $\mathcal{A}$ and $%
\mathcal{B}$, respectively. The maps $\Phi $ then correspond to functors $%
\mathcal{F}$ from $\mathcal{A}$ to $\mathcal{B}$ but since the maps are, in
general, not algebra morphisms, these functors will, in general, not
preserve the quasitensor structure. Since we have a multiplicative structure
on both, $\mathcal{A}$ and $\mathcal{B}$, a suitable class of functors $%
\mathcal{F}$ is endowed with the structure of a bialgebra category in the
sense of \cite{CrFr} (much the same way a suitable class of complex valued
functions on a group is - by inducing the multiplication from the codomain
pointwise - endowed with the structure of a Hopf algebra) where a bialgebra
category is, roughly speaking, a monoidal category with a compatible
functorial version of a coproduct (for the precise definition, see the cited
paper). So, a certain class of ncQFTs will have an algebraic formulation in
terms of bialgebra categories. We will call such ncQFTs ``of bialgebra
category type'' and write bcncQFT for them, for short.

We start our argument on renormalization and bcncQFTs by noting a property
of the moduli spaces (in the sense of formal deformation theory of algebraic
structures, where we will by a moduli space of a structure always mean the
connected component, only) of quasitensor categories defined from two
dimensional conformal field theories.

\bigskip

\begin{lemma}
The moduli space of a quasitensor category $\mathcal{C}$ which is defined
from a two dimensional conformal field theory is always equivalent to the
moduli space generated from the braiding and associativity morphism, alone.
\end{lemma}

\proof%
Since we consider only the connected component, we can consider those parts
of the structure which are descretely parametrized as fixed (e.g. the linear
structures on the homomorphism classes remains fixed; this is completely
analogous to the case of formal deformations of an associative algebra where
only the product is deformed but the underlying linear space - since it is
discretely parametrized by dimension - remains fixed). So, the remaining
structures which can be deformed are the composition, the tensor product,
the braiding, and the associativity morphism for the tensor product.

Now, the moduli space of a two dimensional conformal field theory is locally
parametrized by observables of the theory itself and these, in turn, are in
one to one correspondence to the states of the theory. But then - since
states have to be taken to states by the deformation - every sufficiently
small deformation can be described as a functor on $\mathcal{C}$.
Consequently, the deformations of the composition can always be absorbed
into the other three structures. With the same argument, the deformations of
the tensor product can be assumed to be of the kind of a twist and therefore
can be absorbed into deformations of the associativity morphism. Hence, what
remains are deformations of the braiding and the associativity morphism. 
\endproof%

\bigskip

In the original definition of Drinfeld (see \cite{Dri}), the
Grothendieck-Teichm\"{u}ller group $GT$ is defined from formal deformations
of the braiding and the associativity morphism of a quasitensor category.
Since renormalization is understood as an action on a formal moduli space of
QFTs, the above lemma immediately suggests that the renormalization group
flow induces on the space of two dimensional conformal field theories (part
of) the flow generated by the half group counterpart (see \cite{Dri}) of $GT$
(since by the above lemma there is no other freedom in deforming $\mathcal{C}
$ than the transformations used in Drinfeld's definition). The principle
idea we use, here, is just the treatment of the moduli space of QFTs of
Wilson as a formal moduli space (in the sense of the formal deformation
theory of algebraic structures) of the algebraic formulation of low
dimensional QFTs.

\bigskip

\begin{remark}
The decisive use of the state versus local observable correspondence of two
dimensional conformal field theories shows that we can not necessarily
expect this conclusion on a link between the half group counterpart of $GT$
and the renormalization group flow to hold for the general QFT case.
\end{remark}

\bigskip

In analogy to \cite{Dri}, one of us has introduced in \cite{Sch} from the
possible formal transformations of a braiding and associativity morphism on
a bicategory, plus the two corresponding dual structures for the functorial
coproduct, a noncommutative and noncocommutative self-dual Hopf algebra $%
\mathcal{H}_{GT}$. In complete analogy to the quasitensor category case, one
has the following lemma, then.

\bigskip

\begin{lemma}
The moduli space of a braided and cobraided bialgebra category, arising from
a bcncQFT in the way described at the beginning of this section, is always
equivalent to the moduli space generated from the braiding and associativity
morphism and the corresponding dual structures, alone.
\end{lemma}

\proof%
Completely analogous to the above case. 
\endproof%

\bigskip

In conclusion, one expects that renormalization of a bcncQFT has to be
linked to the Hopf algebra $\mathcal{H}_{GT}$. So, there is evidence from
the general algebraic properties of ncQFTs, too, that renormalization of
such theories should be linked to a noncommutative and noncocommutative Hopf
algebra. The fact that $\mathcal{H}_{GT}$ is self-dual means that it is a
much more noncommutative object than the usual quantum group examples. The
simple deformation of $\mathcal{H}_R$ which we studied in the previous
section, can therefore only expected to be linked to very simple toy models.
For a physically realistic class of ncQFTs we have to expect a much more
complicated Hopf algebra structure.

\bigskip

\section{Conclusion}

We have shown the existence of a noncommutative deformation of the Hopf
algebra of Connes and Kreimer. A left antipode was explicitly constructed
while the existence of a full antipode was only given by an abstract
argument. Surely, the study of deformations of the Connes-Kreimer algebra
and their possible relation to renormalization of quantum field theories on
noncommutative spaces deserves further study. Even if such an approach turns
out to be relevant to renormalization on noncommutative spaces, it is in no
way clear if there exists a canonical noncommutative deformation of the
Connes-Kreimer algebra taking this role or if different deformations
correspond to renormalization schemes for different families of field
theories on noncommutative spaces. The arguments in section 3 show that one
should expect $\mathcal{H}_{GT}$ to be linked to the class of bcncQFTs.

We conclude by one additional remark. The appearance of powers of $q$ given
by a certain index (here: the number of vertices) as the deformation factors
of products and the deformation of the unit element are two features of $%
\mathcal{H}_{R,q}$ which are well known to readers who are acquainted with
trialgebraic deformations of quantum groups (basically, an anewed
deformation quantization of quantum groups, see \cite{GS2000a}, \cite
{GS2000b}). Since trialgebras are linked to the symmetry properties of field
theories on noncommutative spaces (see \cite{GS2000c}), the deformation
theory of $\mathcal{H}_R$ seems nicely to fit in as another element of the
structural properties of such theories.

\bigskip

\textbf{Acknowledgements:}

H.G. was supported under project P11783-PHY of the Fonds zur F\"{o}rderung
der wissenschaftlichen Forschung in \"{O}sterreich. K.G.S. thanks the
Deutsche Forschungsgemeinschaft (DFG) for support by a research grant and
the Erwin Schr\"{o}dinger Institute for Mathematical Physics, Vienna, for
hospitality. Besides this, very special thanks go to Cesar Gomez who posed
the question to us if anybody had ever considered the possibility if
noncommutative deformations of the Connes-Kreimer algebra might be related
to renormalization of quantum field theories on noncommutative spaces and
thereby directly stimulated the present research.

\end{document}